\documentclass[twocolumn]{autart}
\pdfoutput=1
\usepackage{graphicx}
\usepackage{float}
\usepackage{amsfonts}
\usepackage{amsmath}
\usepackage{amssymb,latexsym}
\usepackage{subfig}
\usepackage{natbib}
\newcommand*{\QEDA}{\hfill\ensuremath{\blacksquare}}

\begin{document}

\begin{frontmatter}

\title{Reachability for Partially Observable Discrete Time Stochastic Hybrid Systems \thanksref{footnoteinfo}}
\thanks[footnoteinfo]{This research was funded by National Science Foundation (NSF) Career Award CMMI-1254990, NSF Award CPS-1329878, NSA Science of Security Lablet at North Carolina State University (subaward to the University of New Mexico), and start-up funding from the Department of Electrical and Computer Engineering, University of New Mexico}

\author{Kendra Lesser\corauthref{authorinfo}}\ead{lesser@unm.edu},
\author{Meeko Oishi}\ead{oishi@unm.edu}
\address{Dept. of Electrical and Computer Engineering, University of New Mexico, Albuquerque, USA}
\corauth[authorinfo]{Corresponding Author. Tel.:+1 505 353 2424; Fax: +1 505 277 0299.}

\begin{abstract}

When designing optimal controllers for any system, it is often the case that the true state of the system is unknown to the controller, for example due to noisy measurements or partially observable states.  Incomplete state information must be taken into account in the controller's design in order to preserve its optimality.  The same is true when performing reachability calculations.  To estimate the probability that the state of a stochastic system reaches, or stays within, some set of interest in a given time horizon, it is necessary to find a controller (or at least prove one exists) that drives the system to that set with maximum probability.  This controller, however, does not have access to the true state of the system.  To date, little work has been done on stochastic reachability calculations with partially observable states. What work has been done relies on converting the reachability optimization problem to one with an additive cost function, for which theoretical results are well known.  Our approach is to preserve the multiplicative cost structure when deriving a sufficient statistic that  reduces the problem to one of perfect state information.  Our transformation includes a change of measure that simplifies the distribution of the sufficient statistic conditioned on its previous value.  We develop a dynamic programming recursion for the solution of the equivalent perfect information problem, proving that the recursion is valid, an optimal solution exists, and results in the same solution as to the original problem.  We also show that our results are equivalent to those for the reformulated additive cost problem, and so such a reformulation is not required.

\end{abstract}

\end{frontmatter}

\section{Introduction}

The concept of stochastic reachability enables calculation of the probability that the state of a dynamical system will reach a desired set within a given time horizon.  Alternately, the safety of the system may be considered by examining the probability that the system remains within some safe region.  This problem has gained particular attention in the area of stochastic hybrid systems, which are composed of both continuous and discrete co-evolving states.

The reachabilty problem has been studied mainly for deterministic hybrid systems.  In particular, level set methods have been used to approximate the solution to an appropriate Hamilton-Jacobi-Bellman equation, as in \cite{oishi1}, and \cite{Mitchell2}.  Interest in stochastic systems led to developments in the area of reachability for continuous time hybrid systems whose dynamics update stochastically.  Preliminary results focused on stochastic differential equations combined with deterministic jumps between discrete states in \cite{Hu1}, and was extended to allow for random discrete jumps as well in \cite{Prand1}.  More recently, \cite{Bujor2} examined the continuous time reachability problem as an optimal stopping problem, characterized in terms of variational inequalities.  \cite{Esfahani} also considered the reachability problem as a stochastic optimal control problem with discontinuous payoff functions, and developed a weak dynamic programming principle for its value function.  However, measurability complications that  make solution strategies for the reachability problem more difficult to characterize and derive  have lead to interest in the discrete time equivalent, which circumvents many of those issues (see \cite{Abate1}).

The discrete time stochastic hybrid system (DTSHS) is presented in \cite{Abate1}, and the reachability problem solved as a stochastic optimal control problem using dynamic programming, based on the theory and techniques presented in \cite{bertsekas}.  \cite{summers} extend \cite{Abate1} to the \emph{reach-avoid} problem, where the objective is to avoid an unsafe region while ultimately arriving at a target set, and establishes its solution using the techniques of \cite{bertsekas}.  Finally, \cite{Kamgar1} build upon both of the above works to allow for a disturbance acting in opposition to the controller (as in a two-player dynamical game).

The goal of our presented work is to extend \cite{Abate1}, \cite{summers}, and \cite{Kamgar1} to the case of a partially observable system, where the controller only has access to noisy (possibly incomplete) measurements of the state.  We specifically examine the reachability problem as presented in \cite{Abate1}, which focuses on the safety problem of keeping the state within a known safe region.  Imperfect state information could ultimately lead to suboptimal control inputs as compared to the case in which the true state of the system is known.  In the case of safety verification and reach-avoid set calculations, if the controller is falsely assumed to have perfect knowledge of the state, the reachability probability may be overestimated.  This is certainly undesirable in the context of safety verification.  

There has been extensive work on hybrid estimation (see for instance \cite{hof1}, \cite{Hwang1}, \cite{Kouts1}), but its application to the reachability problem for hybrid systems is limited.  In particular, \cite{Verm1}, \cite{Verm2} examine a continuous time hybrid system subject to continuous control inputs and both continuous and discrete disturbance inputs, where the discrete mode of the system is unknown.  \cite{Verm2} assumes separation between state estimation and control, and reduces the problem of hidden discrete modes to one of perfect state information, by redefining the state.  However, in the case of reachability for hybrid systems, separation between state estimation and control cannot be assumed optimal, and so state estimation cannot be directly applied to the reachability problem. 

Only very recently has there been a surge in work on the reachability problem for partially observable DTSHS (see \cite{Ding2013}).  Although the reachability problem was originally presented in terms of a multiplicative cost function (\cite{Abate1}), \cite{Ding2013} rewrites it as a terminal cost function,  by appending to the state of the hybrid system a binary variable representing whether the state has remained within the desired region up to the previous time step. The partially observed control problem can be reformulated in terms of a sufficient statistic, which  encapsulates and condenses all necessary information for the control of a system.  Thus the problem is recast as one of perfect information, for which solution strategies are well known. 

\cite{Ding2013} make use of the fact that for an additive cost function, the posterior distribution, or probability density of the state given all available information (observations, control inputs) up to the present, provides sufficient information for control of the system (this result is derived in, e.g., \cite{bertsekas}).  Indeed, inspired by this approach, \cite{Tkachev2013} reformulate the reachability problem more generally as an additive cost optimal control problem, although they do not discuss the partially observed case.

While we examine the same problem as \cite{Ding2013}, our derivations preserve the multiplicative cost structure of \cite{Abate1}.   For a nonadditive cost function the posterior distribution is no longer sufficient (see \cite{shiryaev}), and a different sufficient statistic for reducing the problem to one of perfect information must be derived.  However, we will show that while our approach differs from that of \cite{Ding2013}, our results are in fact nearly identical.  Although we preserve the multiplicative cost function, leading to a seemingly more complex problem, the additive cost formulation effectively moves the complexity from the cost function to the modified state of the system.  The posterior distribution of the new state is actually the same as the distribution produced by the sufficient statistic we will derive, so that ultimately the only advantage to the additive cost formulation is its familiar and well-studied form.  Further, we make use of a change of measure in formulating the sufficient statistic that enables easier calculation (by simplifying the distribution of the state of the sufficient statistic conditioned on its previous value). 

To derive a sufficient statistic while preserving the multiplicative cost formulation of the reachability problem, we draw mainly upon theoretical work done in the context of partially observable risk-sensitive stochastic optimal control problems.  The risk-sensitive control problem minimizes the \emph{exponential} of a sum of costs, rather than a sum of costs, so that the cost objective is in fact nonadditive.   In particular, \cite{James} derived a sufficient statistic for a partially observable discrete-time nonlinear system, which was further analyzed and extended in the context of a partially observable Markov decision process (POMDP) in \cite{Gauch}.  In the latter, the state, observation, and control took values from finite, discrete sets, whereas in the former all values were continuous.  As an aside, this highlights how such problems can be regarded often interchangeably as a stochastic optimal control problem or as a Markov decision process (MDP), and results from one field usually carry over to the other, assuming the system dynamics follow the Markov property.  We present here the reachability problem using a control theory framework, and relate the controlled stochastic hybrid system to an MDP.  

Motivated mainly by \cite{James}, we derive a) a sufficient statistic for the multiplicative reachability cost function with hybrid state dynamics, and b) the dynamic programming (DP) equations to solve the reachability problem in terms of the sufficient statistic.  We also introduce a change of measure to the hybrid space, so that the observations are independent of the state of the system (and are in fact independent and identically distributed).  This makes for simpler dynamic programming equations, and should aid in computation and simulations.   The novelty of this work is therefore 1) preservation of a \emph{multiplicative} cost function for the reachability problem to verify the safety of a partially observable DTSHS, 2) introduction of a change of measure to make the observations independent and identically distributed 3) derivation of a sufficient statistic to convert the partially observed problem to a fully observed one, and 4) validation of a DP recursion to solve the reachability problem in terms of the sufficient statistic.   Our main focus is on the theoretical foundations for solving the reachability problem for a partially observable DTSHS.  We note that  implementation of our technique will require further work in approximation strategies as well as in special classes of systems in which exact solutions are available.  

The paper is organized as follows.  First, we review the characteristics of a DTSHS, then extend it to include a hybrid observation space in Section \ref{DTSHS}.  We define the reachability problem, as in \cite{Abate1}, and derive a sufficient statistic, a recursion to update the state of the sufficient statistic, and DP equations for the reachability problem in Section \ref{RCP}.  Here we will also elaborate upon the technique presented in \cite{Ding2013}, and its relation to our own  method.  In Section \ref{example}, we describe two examples of partially observable discrete time stochastic hybrid systems, demonstrate how to reformulate them in terms of our sufficient statistic, and discuss some of the computational challenges as well as possible solution strategies.  Concluding remarks are given in Section \ref{conc}.

\section{Discrete Time Stochastic Hybrid Systems}\label{DTSHS}

A hybrid system is characterized by a set of both discrete and continuous states with interacting dynamics. The discrete state may affect the evolution of the continuous dynamics, and the continuous dynamics may affect when the discrete state changes.  In the case of a discrete time stochastic hybrid system (DTSHS), both the discrete and continuous dynamics may be characterized by stochastic kernels, the product of which determines the stochastic transition kernel governing the combined discrete/continuous state of the system.  We present a slightly modified definition of a DTSHS first introduced in \cite{Abate1}.

\begin{defn}\label{dtshs}(Discrete Time Stochastic Hybrid System). A DTSHS is a tuple $\mathcal{H} = (Q,S,\mathcal{U},T_x,T_q)$ where
	
\begin{enumerate}
\item
$X \subseteq \mathbb{R}^n$ is a set of continuous states
\item
$\mathcal{Q} = \{q_1, q_2, ... q_{N_q}\}$ is a finite set of discrete states with cardinality $N_q$, and $\mathcal{S} = X\times \mathcal{Q}$ is the hybrid state space
\item $\mathcal{U}$ is a compact Borel space which contains all possible control inputs affecting discrete and continuous state transitions
\item$T_x : \mathcal{B}(\mathbb{R}^n)\times \mathcal{S} \times \mathcal{U} \rightarrow [0,1]$ is a Borel-measurable stochastic kernel which assigns a probability measure to $x_{k+1}$ 	
given $s_k=(x_k,q_k),u_k, q_{k+1}\,\forall\, k$: $T_x(dx_{k+1}\in B\mid q_{k+1}, s_k,u_k)$ where $B\in \mathcal{B}(\mathbb{R}^n)$, the Borel $\sigma$-algebra on $\mathbb{R}^n$.
\item $T_q : \mathcal{Q}\times \mathcal{S} \times \mathcal{U} \rightarrow [0,1]$ is a discrete transition kernel assigning a probability distribution to $q_{k+1}$ given $x_k,q_k,u_k,\, \forall \, k$.
\end{enumerate}
\end{defn} 

Kernels $T_x$ and $T_q$ can be combined for ease of notation to produce the hybrid state transition kernel
\begin{equation}\label{statetrans}
\tau(ds' \mid s,u) = T_x(dx' \mid x, q, u, q')T_q(q'\mid x,q,u) \\
\end{equation}
The discrete state $q_{k+1}$ update depends on $q_k$, $x_k$ and $u_k$, and the continuous state $x_{k+1}$ update depends on $x_k$, $u_k$, and according to the specific problem may also be governed by $q_k$, $q_{k+1}$, or both.  For ease of notation we assume that the discrete state updates first, and the updated discrete state affects the continuous state, i.e. that $T_x(dx_{k+1} \mid x_k, u_k, q_{k+1})$, although modifying $T_x$ to include $q_k$ would not alter any subsequent results.  Hence the hybrid system can also be modeled as an MDP with state space $\mathcal{S}$, control space $\mathcal{U}$, and transition function $\tau(s'\mid s, u)$.

\subsection{Partially Observable DTSHS}\label{PODTSHS}

We assume the hybrid process $\left\{x_k, q_k\right\}$ is not measured, but rather that an observation process $y_k = (y_k^x, y_k^q)$ is available.   
The observations $y_k^x \in \mathbb{R}^n$ of the continuous state, and observations $y_k^q\in \mathcal{Q}$ of the discrete state, are assumed independent given the true state $(x_k, q_k)$, so that
\begin{align}
y_k^x &= h(x_k,u_{k-1}) + v_k \label{obs1} \\
y_k^q &\sim Q_{q_k,y^q}(u_{k-1}) \label{obs2}
\end{align}
 The continuous state observation $y_k^x$ is subject to additive noise $v_k$, which is independent and identically distributed with positive density $\phi(v)$ (i.e. Gaussian), and the function $h$ is assumed bounded and continuous, as in \cite{James}.  The distribution of the discrete state observation $y_k^q$ follows the discrete map $ Q_{q,y^q}(u): \mathcal{Q}\times\mathcal{Q}\times \mathcal{U}\rightarrow [0,1]$, so that 
 $P[y_k^q = n \mid q_k=q,\,u_{k-1}=u] = Q_{q,n}(u)$.  The filtrations $\mathcal{G}_k$ and $\mathcal{Y}_k$ are generated by the sequences $\{s_0,\dotsc,s_{k},y_1,\dotsc,y_{k-1}\}$ and $\{y_1,\dotsc,y_k\}$, respectively.  Denote $i_k = (y_1,\dotsc,y_k, u_0,\dotsc,u_{k-1})\in \mathcal{Y}_k\times\mathcal{U}_k = \mathcal{I}_k$, with $\mathcal{U}_k$ the $k$-times product space of $\mathcal{U}$, as the vector of information available at time $k$.   The information vector $i_k$ is used to make the control decisions $u_k$ through a control policy $\pi= (\mu_0,\dotsc,\mu_{N-1})$, where $\mu_k: \mathcal{I}_k \rightarrow \mathcal{U}$ is a function mapping the space of available information, $\mathcal{I}_k$, into the control space $\mathcal{U}$ for all $k=0,\dotsc,N-1$. 

We also assume an initial Borel-measurable density on $s_0=(x_0,q_0)$, $s_0\sim\rho(x,q) \in P(\mathcal{S})$, i.e. that $\rho$ lies in the space of all probability measures on $\mathcal{S}$.  Finally, based on $\rho$, $\tau$, $\phi$, and $Q(u)$, we obtain a probability measure $\mathbb{P}^{\pi}$ induced by the control policy $\pi$ defined over the full state space $\Omega$, which includes $s_k$ and $y_k$ for all $k$.
We can therefore model the hybrid system with observations as a POMDP.


\section{Reachability Control Problem}\label{RCP}

\subsection{Cost Function}

We wish to find a control policy that maximizes the probability that the true state of the system stays within some safe or desired set for a finite time horizon.  As in \cite{Abate1}, this problem can be formulated as a stochastic optimal control problem.  For a Borel set $K$, terminal time $N$, and predefined policy $\pi$, define the cost function as
\begin{equation}\label{probRA}
r_{K}(\pi) = \mathbb{P}^{\pi}[s_i \in K \, \forall \,i=0,\dotsc,N]
\end{equation}
Since for a random variable $X$, $\mathbb{P}[x \in A] = \mathbb{E}[{\bf 1}_A(x)]$, with $\mathbb{E}$ denoting expected value and indicator function ${\bf1}_A(x)=1$ if $x \in A$ and ${\bf1}_A(x)=0$ otherwise, \eqref{probRA} is rewritten as in \cite{Abate1}:
\begin{equation} \label{ERA}
r_K(\pi) ={\mathbb E}^{\pi} \left[ \prod_{i=0}^{N}{\bf 1}_K(s_i) \right]
\end{equation}
The expected value is taken with respect to the measure $\mathbb{P}^{\pi}$, hence the notation $\mathbb{E}^{\pi}$.  We want to maximize $r_K(\pi)$ with respect to the control policy $\pi$.  The set $\Pi$ of admissible policies will be restricted to non-randomized policies, i.e. in which $\mu_k(i_k)$ generates one control input $u_k$ with probability $1$.   The optimal policy $\pi^*$ is then given by
\begin{equation}\label{optpol}
\pi^* = \arg\sup_{\pi\in\Pi} \left\{r_K(\pi)\right\}
\end{equation}
 We can now formally define the partially observable reachability problem we wish to solve.
\begin{defn}\label{ProblemStatement}
\emph{(Problem Statement)}
Consider a DTSHS (defined in Definition \ref{dtshs}) with observations \eqref{obs1} - \eqref{obs2} and initial distribution $\rho(x,q)\in P(\mathcal{S})$.  Given a safe set $K$ and time horizon $N$ we would like to
\begin{enumerate}
\item
Find the maximal probability of remaining within $K$ for $N$ time steps, given by $\max_{\pi} r_K(\pi)$.
\item
Find the optimal policy $\pi^*$ such that  $\max_{\pi} r_K(\pi) = r_K(\pi^*)$.
\end{enumerate}
\end{defn}

In the case of perfect state information, where the control $u_k$ is a function of $s_k$ rather than $i_k$, Definition \ref{ProblemStatement} can be solved via dynamic programming as presented in \cite{Abate1}.  When only an observation process is available, however, the standard approach is to reformulate the problem as one with perfect information by redefining the state of the system in terms of a \emph{sufficient statistic} (see, e.g. \cite{bertsekas}, \cite{Aoki}) and then solving the equivalent problem using dynamic programming. 

The difficulty in solving Definition \ref{ProblemStatement} is twofold.  First, since the cost function is multiplicative, standard sufficient statistics are not valid (i.e. the sufficient statistic cannot be the posterior distribution of the state at time $k$ given all available information up to time $k$).  Second, the hybrid nature of the dynamics complicates the probability space our problem is defined on.  A new sufficient statistic must be derived, and its corresponding theoretical results carefully extended.

\subsection{Sufficient Statistic}


We will first formally define a sufficient statistic in relation to the multiplicative optimal control problem of Definition \ref{ProblemStatement}, which is modified from Definition 10.6 in \cite{bertsekas}.

\begin{defn}\label{suffStatDef}
A statistic for Definition \ref{ProblemStatement} is a sequence of Borel-measurable functions $(\eta_0, \eta_1,\dotsc,\eta_N)$ with $\eta_k: P(\mathcal{S})\times \mathcal{I}_k \rightarrow \Sigma_k$ where $\Sigma_k$ is a nonempty Borel space, for all $k = 0,\dotsc,N$.  The sequence $(\eta_0,\dotsc\eta_N)$ is a statistic \emph{sufficient for control} if
\begin{enumerate}
\item
There exist Borel-measurable stochastic kernels $\hat{\tau}_k(d\sigma_{k+1}\mid \sigma_k, u_k)$ on $\Sigma_{k+1}$ given $\Sigma_k$, $\mathcal{U}$ such that
\begin{align*}
\mathbb{P}^{\pi}&[\eta_{k+1}(\rho, i_{k+1})\in \underline{\Sigma}_{k+1}\mid\eta_k(\rho, i_k) =\sigma_k, u_k = \overline{u}_k]\\
&\hspace{20 mm} = \hat{\tau}_k(\underline{\Sigma}_{k+1}\mid \sigma_k,  \overline{u}_k )
\end{align*}
for $\mathbb{P}^{\pi}$ almost every $(\sigma_k, \overline{u}_k)$ (i.e. up to a set of measure zero with respect to measure $\mathbb{P}^{\pi}$).
\item
There exist functions $g_k:\Sigma_k\rightarrow [0,\infty)$  such that for all $\rho\in P(\mathcal{S})$, $k=1,\dotsc,N$, and $\pi\in\Pi$ 
\begin{equation*}
 \mathbb{E}^{\pi}\left[\left.\prod_{i=1}^k{\bf 1}_K(s_i)\right| \eta_k(\rho,i_k) =\sigma_k\right]
= g_k(\sigma_k)
\end{equation*}
for $\mathbb{P}^{\pi}$  almost every $\sigma_k$.
\end{enumerate} 
\end{defn}

In other words, the distribution of $\sigma_k = \eta_k(\rho,i_k)$ (a specific value of the sufficient statistic which we refer to as the \emph{information state}) must follow the Markov property, and therefore be updated recursively according to $\sigma_{k-1}$ and $u_{k-1}$.  There also must exist an equivalent cost function whose argument is $\sigma_k$, so that the cost corresponding to a specific policy can be determined solely through the distribution of the information state.   The problem presented in Definition 2 can then be redefined according to the information state $\sigma$, which itself is defined according to the sufficient statistic $\eta$.

We now propose a sufficient statistic for the partially observable reachability problem, and demonstrate that it obeys Definition \ref{suffStatDef}.  First, we introduce the concept of a change of measure, which is used in the derivation of our sufficient statistic to facilitate the analysis and subsequent computation.  The ability to change probability measures stems from the Radon-Nikodym Theorem.
\begin{defn}
The Radon-Nikodym Theorem (see \cite{SteinShak}) states that given two $\sigma$-finite measures $\nu$ and $\mu$ on a measurable space $(\Omega, \mathcal{M})$, if $\mu$ and $\nu$ are absolutely continuous, then there exists a $\mu$-integrable function $f$ on $\Omega$ such that
\begin{equation*}
\nu(E) = \int_E f(\omega)\,d\mu(\omega)
\end{equation*}
The function $f$ is referred to as the Radon-Nikodym derivative, and is written as $\frac{d\nu}{d\mu}$.
\end{defn}

Essentially, for two probability measures $\nu$ and $\mu$, defined on the same space $(\Omega, \mathcal{M})$ and that satisfy $\nu(E) = 0$ whenever $\mu(E)=0$ for all $E\in\mathcal{M}$, we know that
\begin{equation*}
\mathbb{E}^{\nu}[h(\omega)] = \mathbb{E}^{\mu}[f(\omega)h(\omega)]
\end{equation*}
for any $\mathcal{M}$-measurable function  $h$.  We can define a change of measure $\mathbb{P}^{\dag}$ from the existing measure $\mathbb{P}^{\pi}$ on our space $\Omega$, with $\mathcal{M}$ being the Borel $\sigma$-algebra on $\Omega$, so long as the continuous observation process is nowhere zero, and the discrete observation is nowhere zero on $\mathcal{Q}\times\mathcal{Q}\times\mathcal{U}$ (which would occur if certain discrete states were perfectly observable).  Following  \cite{James} and \cite{Elliot93}, we define the Radon-Nikodym derivative $\Lambda_k$ as 
\begin{equation}\label{RNder}
\left. \frac{d\mathbb{P}^{\pi}}{d\mathbb{P}^{\dag}}  \right|_{\mathcal{G}_k} = \Lambda_k
\end{equation}
where
\begin{equation*}
\Lambda_k =  \prod_{l=1}^k \frac{\phi(y_l^x - h(x_l,u_{l-1}))Q_{q_l,y_l^q}(u_{l-1})}{\phi(y_l^x)\frac{1}{N_q}}
\end{equation*}
 However, in contrast to \cite{James},  we must contend with two separate observation processes, one continuous and one discrete.  Note that in \eqref{RNder} we restrict the derivative to the filtration $\mathcal{G}_k$ rather than the full state space $\Omega$, which allows us to update the derivative as the hybrid process evolves.
\begin{lem}\label{iidmeas}
Under $\mathbb{P}^{\dag}$, the processes $\{y_k^x\}$ are independent and identically distributed (i.i.d.) with density $\phi$, and the processes $\{y_k^q\}$ are i.i.d. with uniform density $\frac{1}{N_q}$.
\end{lem}

The proof of Lemma \ref{iidmeas} is in the Appendix.  We claim that under the measure $\mathbb{P}^{\dag}$ we can define a sufficient statistic for the reachability control problem.
\begin{equation}\label{sigma}
 \eta_k(\rho, i_k)=  \mathbb{E}^{\dag}\left[\left. {\bf 1}_q(q_k){\bf 1}_{x}(x_k) \prod_{i=1}^{k-1}{\bf 1}_K(s_i) \Lambda_k \right| i_k\right]
\end{equation}
The sufficient statistic $(\eta_0, \eta_1,\dotsc,\eta_N)$ generates a sequece of unnormalized probability densities on the state $s_k$, condtioned on the information vector $i_k$, and where the dependence on $\rho$ is implicitly defined in the measure $\mathbb{P}^{\dag}$.  The information state $\sigma_k$ is therefore a modification of the posterior distribution, and represents an unnormalized conditional density of the current state joined with the probability that all previous states are in $K$, given a specific $i_k$.  In order to show that our sufficient statistic \eqref{sigma} satisfies conditions (1) and (2) of Definition \ref{suffStatDef}, we refer explicitly to the informations state $\sigma_k$ rather than the sufficient statistic $\eta_k$, noting that although the sufficient statistic is a function of $\rho$ and $i_k$, the information state itself is a function of the state $s_k$, and its dependence on $\rho$ and $i_k$ is not explicitly stated.  We first show that $\sigma_k$ can be defined recursively via a bounded linear operator $T_{u,y}[\sigma]$, and therefore satisfies (1) of Definition \ref{suffStatDef}.  The proof can be found in the Appendix.

\begin{lem}\label{linop}
There exists a bounded linear operator $T: L^1(\mathcal{S})\rightarrow L^1(\mathcal{S})$ such that $\sigma$ is defined recursively as:
\begin{equation}\label{sigma_rec}
\begin{cases}
\sigma_0 = \rho \\
\sigma_k = T_{u_{k-1},y_k}[\sigma_{k-1}]
\end{cases}
\end{equation}
where $T_{u,y}[\sigma]$ is given by
\begin{align}
T_{u,y}[\sigma] &= \sum_{q^-\in\mathcal{Q}}N_qQ_{q,y^q}(u)\int_{\mathbb{R}^n} {\bf 1}_K(x^-,q^-) \notag \\
&\times  \frac{\phi(y^x - h(x,u))}{\phi(y^x)} \tau(x,q\mid x^-,q^-,u)\sigma(x^-,q^-)\,dx^-  \label{T}
\end{align}
In addition, $\sigma_k\in L^1(S)$ for all $k$ since $\sigma_0 = \rho \in L^1(\mathcal{S})$ and $T$ maps $L^1$ into $L^1$.
\end{lem}

 The stochastic kernel $\hat{\tau}_k$ for the distribution of $\sigma_{k+1}$ given $\sigma_k$ and $u_k$ can be written in terms of the new measure $\mathbb{P}^{\dag}$:
\begin{equation}\label{tauhat}
\hat{\tau}_{k}(\underline{\Sigma}_{k+1}\mid \overline{\sigma}_k, \overline{u}_k) =\sum_{y^q\in\underline{Y}^q}\int_{\underline{Y}^x} \frac{1}{N_q}\phi(y^x)\,dy^x
\end{equation}
with $\underline{Y}^q\times\underline{Y}^x = \{(y^q, y^x) : T_{\overline{u}_k, y^x, y^q}[\overline{\sigma}_k] \in \underline{\Sigma}_{k+1}\}$.  

Next, we rewrite the cost function \eqref{ERA} in terms of the information state $\sigma$, for part (2) of Definition \ref{suffStatDef}.  Since the indicator function ${\bf 1}_K(s)$ is in the space $L^{\infty}(\mathcal{S})$, the inner product of $\sigma$ and ${\bf 1}_K$ is a well defined bounded linear operator, given by
\begin{equation*}
\langle\sigma, {\bf 1}_K\rangle = \sum_{q\in \mathcal{Q}}\int_{\mathbb{R}^n} \sigma(x,q) {\bf 1}_K(x,q) \, dx
\end{equation*}

The functions $g_k$ in Definition \ref{suffStatDef} can be defined as 
\begin{align}
g_k(\sigma_k)=&\mathbb{E}^{\pi}\left[\left.\prod_{i=1}^k {\bf 1}_K(s_i) \right| \eta_k(\rho, i_k) = \sigma_k \right] \notag\\
=&\mathbb{E}^{\dag}\left[ \left. {\bf 1}_K(s_k) \prod_{i=1}^{k-1}{\bf 1}_K(s_i) \Lambda_k \right|i_k\right] \notag\\
 =& \langle \sigma_k,  {\bf 1}_K \rangle \label{gfunc}
\end{align}

\subsection{Equivalence to Perfect State Information Problem}
Under the measure $\mathbb{P}^{\dag}$, we can rewrite the cost function \eqref{ERA} as $r_K(\pi)=\mathbb{E}^{\dag}\left[  \Lambda_N \prod_{i=0}^{N}{\bf 1}_K(x_i) \right]$.  We define 	 
\begin{equation}\label{sigcost}	 
\overline{r}_K({\overline{\pi}})=  \mathbb{E}^{\dag}\left[ \langle \sigma_N,{\bf 1}_K \rangle \right]
\end{equation}
as the equivalent cost function for the reachability problem in terms of $\sigma$, using $\overline{\pi}$ to denote a policy in terms of $\sigma$ (whereas $\pi$ denotes a policy for the partially observed case in terms of $i_k$).   Note that for a fixed vector of control inputs (i.e. open loop), ${\bf u}=[u_0,\dotsc,u_{N-1}]$, $r_K({\bf u}) = \overline{r}_K({\bf u})$:
\begin{align*}
r_K({\bf u}) &= \mathbb{E}^{\dag}\left[  \Lambda_N \prod_{i=0}^{N}{\bf 1}_K(s_i) \right] \\
&=  \mathbb{E}^{\dag}\left[ \mathbb{E}^{\dag}\left[ \left. \Lambda_N \prod_{i=0}^{N}{\bf 1}_K(s_i) \right| i_N\right]\right] \\
&= \mathbb{E}^{\dag}\left[ \langle \sigma_N,{\bf 1}_K \rangle \right] \\
&= \overline{r}_K({\bf u})
\end{align*}
where the third line follows from \eqref{gfunc}.  A recursion for $\overline{r}_K(\overline{\pi})$ is given as in \cite{James} by
\begin{equation}\label{DP}
\begin{cases}
V_N^{\overline{\pi}}(\sigma) = \langle \sigma, {\bf 1}_K \rangle \\
V_k^{\overline{\pi}}(\sigma) =  \mathbb{E}^{\dag}\left[ V_{k+1}^{\overline{\pi}}(T_{\overline{\mu}_k(\sigma),y_{k+1}}[\sigma])\right]
\end{cases}
\end{equation}
where
\begin{align}
&\mathbb{E}^{\dag}\left[ V_{k+1}^{\overline{\pi}}(T_{\overline{\mu}_k(\sigma),y_{k+1}}[\sigma])\right] \notag \\
&\hspace{5 mm}= \int_{\Sigma} V_{k+1}(\sigma')\hat{\tau}(\sigma'\mid \sigma, \overline{\mu}_k(\sigma)) \,d\sigma'
 \label{expectedV} 
\end{align}
\begin{align}
&\hspace{5 mm} = \sum_{\mathcal{Q}} \int_{\mathbb{R}^n} V_{k+1}(T_{\overline{\mu}_k(\sigma),y_{k+1}}[\sigma])\frac{1}{N_q}\phi(y_{k+1}^x)\,dy^x \label{expectedV2}
\end{align}

so that $V_0^{\overline{\pi}}(\rho) = \overline{r}_K(\overline{\pi})$.  Next, we provide two theorems: 1) a dynamic programming algorithm to find the optimal solution to $\sup_{\overline{\pi}}\overline{r}_K(\overline{\pi})$, and the optimal policy $\overline{\pi}^*=\arg\sup_{\overline{\pi}}\overline{r}_K(\overline{\pi})$ as a function of the information state $\sigma$, and 2) proof that this optimal policy has the same value as the optimal policy for the partially observed case.  The proofs of Theorems \ref{DPthm} and \ref{equival} are provided in the Appendix. 
\begin{thm}\label{DPthm}
Using the recursion \eqref{DP}, the dynamic programming equations
\begin{equation}\label{optDP}
\begin{cases}
V_N^*(\sigma) = \langle \sigma, {\bf 1}_K \rangle \\
V_k^*(\sigma) = \sup_{u\in\mathcal{U}} \mathbb{E}^{\dag}\left[ V_{k+1}^*(T_{u,y_{k+1}}[\sigma])\right]
\end{cases}
\end{equation}
produce $V_0^*(\rho) = sup_{\overline{\pi}\in\overline{\Pi}}\overline{r}_K(\overline{\pi})$, where $V_k:\Sigma_k\rightarrow [0,\infty)$. For $\sigma$ normalized, we have  $V_k:\Sigma_k\rightarrow [0,1]$.  Furthermore, setting 
\begin{equation}\label{DPpol}
\overline{\mu}_k^*(\sigma) = \arg\sup_{u\in \mathcal{U}}  \mathbb{E}^{\dag}\left[ V_{k+1}^*(T_{u,y_{k+1}}[\sigma])\right]
\end{equation}
for all $k= 0,\dotsc,N-1$ gives the optimal policy $\overline{\pi}^* = (\overline{\mu}_0^*,\overline{\mu}_1^*,\dotsc,\overline{\mu}_{N-1}^*)$, where $\overline{\mu}_k: \Sigma_k\rightarrow \mathcal{U}$.
\end{thm}

\begin{thm}\label{equival}
If $u_k^* = \overline{\mu}_k^*(\sigma_k)$ is optimal as defined in Theorem \ref{DPthm}, then $u_k^*$ is also optimal for the partially observable problem of Definition \ref{ProblemStatement}, and can be written as $u_k^* =\mu_k^*(i_k) = \overline{\mu}_k^*(\eta(\rho, i_k)) = \overline{\mu}_k^*(\sigma_k)$.  For $\pi^* = (\mu_0^*,\mu_1^*,\dotsc,\mu_{N-1}^*)$, $\overline{r}_K(\overline{\pi}^*) = \sup_{\pi}r_K(\pi) = r_K(\pi^*)$.
\end{thm}

These results guarantee that we can solve \eqref{DPpol} as a fully observed problem for each $k$ in terms of the new state $\sigma$, and generate a policy $\overline{\pi}$ in which the optimal action at each time step $k$ is given as a function only of the information state at time $k$.  Calculating the optimal policy $\overline{\pi}^*$ and optimal value $\overline{r}_K(\overline{\pi}^*)$ gives us the optimal policy $\pi^*$ and optimal value $r_K(\pi^*)$. 

\subsection{Relationship to Additive Cost Formulation}

The sufficient statistic \eqref{sigma} modifies the posterior distribution to include the probability of all previous states being in the set $K$.  Were the sufficient statistic \eqref{sigma} derived without the change of measure \eqref{RNder}, it would be identical (aside from a normalizing constant) to the sufficient statistic for the additive cost function formulation in \cite{Ding2013} (see equation (14)).  In \cite{Ding2013}, by extending the state to include a binary variable that represents whether or not the system has remained within $K$ up to the previous time step, the posterior distribution  is also the distribution of the current state $s_k$, coupled with the distribution of all previous states being in $K$.  The transition kernel for the modified state in \cite{Ding2013}, equation (5),  incorporates an indicator function that signals whether the state remained within the safe set at the previous time step.  The prediction and update steps for a Bayesian filter (see equations (11) - (13) of \cite{Ding2013}) are used to express the sufficient statistic (14), which, aside from the change of measure and normalization, is the same as \eqref{sigma}.

Next, the terminal payoff (15) of \cite{Ding2013} expresses the probability that the final state is within set $K$, and that all previous states are within set $K$, given the probability distribution of all previous states being in $K$ as well as the current distribution of the final state.  Were the terminal payoff written in terms of the original state, it would be identical to \eqref{gfunc} for $k=N$.  Prop. 3 in \cite{Ding2013}, which describes the solution of the terminal payoff, iteratively evaluates the expected value as in Theorem \ref{DPthm} here (although integrating the expected value over $\hat{\tau}_k(\sigma_{k+1}\mid \sigma_k, u_k)$ does not reduce to \eqref{expectedV2} as ours does).  Thus, formulating the cost function as either multiplicative or additive ultimately does not alter the end result.

\section{Case Studies and Computational Issues}\label{example}

We provide two examples of partially observable hybrid systems to demonstrate the use of the sufficient statistic in their solution, then discuss computational challenges.  Since solving \eqref{optDP} requires looping over all functions $\sigma\in L^1(\mathcal{S})$, an infinite space, we can only hope to use \eqref{optDP} as a practical solution method for special cases in which $\sigma$ can be defined over a finite subspace of $L^1$.

\subsection{Temperature Regulation}

A stochastic version of the benchmark temperature regulation problem with perfect state information is presented in \cite{Abatecomp}.  We consider the case of one heater, which can either be turned off, or turned on to heat one of $n$ rooms.  The average temperature of room $i$ at time $k$ is given by the continuous variable $x_i(k)$, and the discrete state $q(k)=i$ indicates room $i \in \{1,\dotsc, n\}$ is heated at time $k$, and $q(k)=0$ denotes the heater is off.  The stochastic difference equation governing the average temperature for room $i$ is given by
\begin{align*}
x_i(k+1) = (1-b_i)x_i(k)+\sum_{i\not= j}&a_{i,j}(x_j(k)-x_i(k))\\
&+c_i h_i+b_ix_a + v_i(k)
\end{align*}
with constants $a_{i,j}$, $b_i$, $c_i$, and $x_a$, $v_i(k)$ that are i.i.d. normally distributed random variables with mean zero and variance $v^2$, and $h_i=1$ for $q(k)=i$ and $h_i=0$ otherwise.  The control input is given by $u(k)\in\mathcal{U}$ with $\mathcal{U}=\{0,1,\dotsc,n\}$, but the chosen control is not always implemented with probability $1$.  Instead, $q(k)$ is updated probabilistically, dependent on $u(k-1)$ and $q(k-1)$, with transition function $T_q(q(k+1)\mid q(k), u(k))$.  So while function $\overline{\mu}_k(\sigma_k)$ deterministically returns a single control input, control input $u(k) = \overline{\mu}_k(\sigma_k)$ may not always be implemented.

The exact average temperature in each room is unknown, and only a noisy measurement of each room's temperature is available to the controller.  The controller does, however, know which room is being heated at time $k$ (i.e. $q(k)$ is perfectly observed).  Then the observation $y(k) = (y^x(k), y^q(k))$ with $y^x(k) = [y_1^x(k)\dotsc y_n^x(k)]^T$ is given by
\begin{align*}
y^x_i(k) &= x_i(k) + w_i(k) \\
y^q(k) &= q(k)
\end{align*}
with $w_i(k)$ i.i.d. normally distributed with mean zero and variance $w^2$ (so that the distribution $\phi(w)$ is Gaussian).  The transition matrix $Q(u)$ is the identity matrix for all $u$, so that $Q_{q,y^q}(u) = {\bf 1}_q(y^q)$. Because the discrete state is perfectly observed, we do not keep track of a discrete observation, and it is not included in the sufficient statistic.

It is desirable to keep the temperature of each room between $17.5$ and $22$ degrees celsius at all times, producing the safe region $K=[17.5,22]\times\dotsc\times[17.5,22]$, which does not depend on the discrete state $q(k)$ (so ${\bf 1}_K(s) = {\bf 1}_K(x)$).  To find the maximum probability that each room stays within the desired temperature range given that the controller only has access to the observations $y(k)$ we reformulate the problem in terms of the information state $\sigma_k$.  We then use the dynamic programming equations \eqref{optDP} given in Theorem \ref{DPthm}, so that
\begin{align*}
V_N^*(\sigma) &= \sum_{q=0}^n\int_{\mathbb{R}^n}{\bf 1}_K(x) \sigma(x,q)\,dx \\  
V_k^*(\sigma) &= \sup_{u\in\mathcal{U}} \int_{\mathbb{R}^n} V_{k+1}^*(T_{u,y}[\sigma])\phi(y^x)\,dy^x
\end{align*}

With $x_i(0)\sim \mathcal{N}(\mu_i,s^2)$ for each $i=1,\dotsc,n$ and $q(0)=0$, then 
\begin{equation*}
\sigma_0(x,q) = {\bf 1}_0(q) \prod_{i=1}^n \rho_i(x_i)
\end{equation*}
for $\rho_i(x)$ Gaussian with mean $\mu_i$ and variance $s^2$.  

However, even for the trivial case of $n=1$ (e.g. a one room system), updating $\sigma_k$ becomes complicated very quickly.  Using Lemma \ref{linop}, we obtain
\begin{align*}
\sigma_1(x,q) &=\frac{\phi(y_1^x-x)}{\phi(y_1^x)}T_q(q\mid q(0)=0,u(0))\\
&\hspace{5 mm}\times \int_KT_x(x\mid x(0), q, u(0))\sigma_0(x(0),q(0)) \,dx(0)
\end{align*}
The main difficulty with solving for $\sigma_1$ is the fact that the integral is evaluated over $K$, as opposed to over $\mathbb{R}$.  
Because of the bounds on the integral, we cannot claim $\sigma_k$ is Gaussian given that $\sigma_{k-1}$ is.  However, because the expression does quite closely resemble a Gaussian distribution, it may be possible to approximate $\sigma_k$ by an un-normalized Gaussian distribution without losing significant accuracy.  We intend to explore this possibility in future works.  Further, we note that there may be classes of systems for which such straightforward sufficient statistics may be found.

%

\subsection{Skid-Steered Vehicle}

A skid-steered vehicle (SSV), modeled as a switched system, is presented in \cite{caldwell}.  The SSV moves according to lateral sticking and sliding of its four wheels.  \cite{caldwell} identify four modes associated with the vehicle: In mode 1, front and rear wheels stick laterally; in mode 2, front wheels stick and rear wheels skid laterally; in mode 3, front wheels skid and rear wheels stick; in mode 4, both front and rear wheels skid laterally.  For each mode, the vehicle's continuous states $X$, $Y$, and $\theta$ are governed by a different set of second order ordinary differential equations (ODEs).  The states $X$ and $Y$ represent the cartesian coordinates for the vehicle's center of geometry, and $\theta$ gives the heading of the vehicle.  We can represent the continuous state of the system by $x = (X,\dot{X},Y,\dot{Y},\theta,\dot{\theta})$, such that $ \dot{x} = f_q(x) = f(x,q)$, with discrete state  $q \in \mathcal{Q}=\{1,2,3,4\}$.  See \cite{caldwell} for the actual expressions for $f_q$, which are too lengthy to reproduce here.  We discretize the ODEs using an Euler approximation method to produce an equivalent discrete time system.  

The control input can be expressed as a command informing the vehicle of what mode it should be in.  If the vehicle responded perfectly, we would have $q_k = u_k$ for the mode at time $k$.  Instead, let us assume that the mode changes behave similarly to the temperature regulation problem above, where the control command is implemented with a certain probability, dependent on the current mode: $T_q(q_{k+1}\mid q_{k}, u_{k+1})$.  The continuous state is assumed to be deterministic given the mode, so that $T_x(x_{k+1}\mid x_k, q_{k},u_k) = {\bf 1}_{f(x_k,q_{k})}(x_{k+1})$.  Finally, assume we have a noisy observation of the continuous state, and have an observation of the mode which is not completely reliable:
\begin{align*}
y_k^x &= x_k + w_k \\
y_k^q &\sim Q_{q,y^q}
\end{align*}
The vector $w_k \in \mathbb{R}^6$ is an i.i.d. sequence of multivariate Gaussians with $w_k\sim\mathcal{N}(0, \mathcal{W})$.   The matrix $Q_{q,y^q}$ is given by
\begin{equation*}
Q_{q,y^q} = \begin{bmatrix} .9& .033 &.033 & .033 \\ .033 &.9&.033&.033\\ .033&.033&.9&.033\\.033&.033&.033&.9\end{bmatrix}
\end{equation*}
Each row corresponds to a given value of $q$, and the probability that $y^q$ takes on each value one to four.  Thus, the probability that the observed mode equals the true mode is $0.9$, and if it is not the true mode, is equally likely to be any of the other three modes.

The safe region $K$ can be defined as a path we would like the vehicle to stay on, which will be defined in terms of bounds on $X$ and $Y$.  For instance, we could define $K$ as a rectangular strip $K=\{X,Y: -1\leq X \leq 1, \, -10\leq Y \leq 10\}$.  Assuming the initial position of the vehicle, $x_0$, is known and equal to $\hat{x}_0\in K$, and the initial mode is independent of $x_0$, uniformly distributed, and represented by $\rho(q_0)$, $\sigma_0$ is given by
\begin{equation*}
\sigma_0(x,q) = {\bf 1}_{\hat{x}_0}(x)\rho(q) = \frac{1}{4}  {\bf 1}_{\hat{x}_0}(x)
\end{equation*}

In this case $\sigma_1$ is easily calculated:
\begin{align*}
\sigma_1(x,q) &= \sum_{q_0=1}^4\int_K 4Q_{q,y_1^q}T_q(q\mid q_0, u_0)\frac{\phi(y_1^x -x)}{\phi(y_1^x)}\\
&\hspace{29 mm} \times {\bf 1}_{f(x_0,q_0)}(x) \sigma_0(x_0,q_0)\,dx_0 \\
&=\sum_{q_0=1}^4Q_{q,y_1^q}T_q(q\mid q_0, u_0)\frac{\phi(y_1^x - f(\hat{x}_0,q_0))}{\phi(y_1^x)}\\
&\hspace{50 mm} \times {\bf 1}_{f(\hat{x}_0,q_0)}(x)
\end{align*}
Thus there are four possible values of $x$ for which $\sigma_1(x,q)$ is nonzero.  Similarly, given an $x_k$ value, $\sigma_{k+1}(x,q)$ will only be nonzero for four values of $x$.
\begin{align*}
\sigma_{k+1}(x,q) &= \sum_{q_{k}=1}^4\int_K 4Q_{q,y_{k+1}^q}T_q(q\mid q_k, u_k)\\
\times &\frac{\phi(y_{k+1}^x -f(x_k,q_k))}{\phi(y_{k+1}^x)}{\bf 1}_{f(x_k,q_k)}(x)\sigma_k(x_k,q_k)\,dx_k
\end{align*}

Even when $\sigma_k$ takes the above seemingly simple form, there is no immediately obvious way to avoid evaluating the value functions for all $\sigma\in L^1$ in order to solve \eqref{optDP}. 

\subsection{Computational Challenges}

Because of the complexity of the hybrid dynamics and cost function, the sufficient statistic and DP equations are computationally intensive.  The DP equations require looping over an infinite state space. No computational work has yet been done on the reachability problem for partially observable DTSHS, and the applicability of current computational work on general DTSHS seems limited.  For instance, discretization procedures for continuous state processes, like those presented in \cite{SoudjaniSiam13} are not immediately applicable because it is much more difficult to grid a continuous function of a continuous state (like our sufficient statistic).  Other methods involve reformulating the reachability problem using chance-constrained optimization, so that the safety constraint is enforced with some probability.  These chance constraints are often evaluated using sampling-based methods (see, e.g., \cite{Vrakop2013}).  Unfortunately, such methods run into the same problems as DP (curse of dimensionality) if formulated as multistage stochastic programs, and it is again not obvious how to extend such methods to the partialy observed case, where dimensionality is an even greater issue.  

For the two examples presented here, one major challenge is circumventing the evaluation of the value functions for all $\sigma\in L^1(\mathcal{S})$ to solve \eqref{optDP}.  One possible alternative is using approximate DP to estimate the value functions $V_k$ by sampling from $y_k$ for each $k$ to get sample trajectories of the $\sigma_k$.  Since via our change of measure the $y_k$ are i.i.d, such sampling should be straightforward.  Each $y_k^x$ is sampled from $\phi(\cdot)$, and each $y_k^q$ is sampled from the uniform distribution on $\{1,\dotsc,N_q\}$.  Some work has explored approximate dynamic programming for DTSHS (see \cite{Kariot2013}), in which the value function is approximated using a linear combination of basis functions, and constraints on the value function are evaluated by sampling from the state space.  It is possible a similar approach could be applied to the partially observed case, where we must sample from the observation space to obtain instances of $\sigma$.

In addition, some work has been done on approximating continuous state POMDPs using point-based value iteration (see \cite{Porta}), albeit in the context of additive cost functions with the belief state as a sufficient statistic.  The method exploits the structure of the value function, and uses Monte Carlo methods to generate a set of samples from the belief space, in order to approximate the value function at a given starting belief state.  Further, \cite{Roy} have applied this to a system with hybrid dynamics.  More recently, there has been a greater focus on solving continuous state POMDPs through approximation and sampling, including a Monte Carlo technique that samples both from the belief space and the state space (\cite{Bai2011}).  Although our cost function and sufficient statistic are different, we are currently working to extend these methods to solve the reachability problem.

 \section{Conclusion}\label{conc}

We have presented a statistic sufficient for the control of a partially observable discrete time stochastic hybrid system, when the objective is to maximize the probability of remaining within a safe set for some finite time horizon.  By redefining the partially observed optimal control problem as one that is fully observed, with state variable $\sigma$ (the information state generate by the sufficient statistic), we are able to define an optimal control policy as a function of $\sigma$.  This control policy is equivalent to the policy defined as a function of the information vector, and leads to the same maximal safety probability.  Further, we showed the equivalence between our approach and one that uses an additive cost function.

The major disadvantage of the sufficient statistic is that the dynamic programming equations must be solved for every possible $\sigma \in L^1$ at every time step.  As a direct solution method, it is seemingly impractical.  However, there may be cases where $\sigma$ can be limited to a subset of $L^1$ so that the dynamic programming equations can be solved.  Further, our choice of measure in defining the sufficient statistic may lend itself well to approximate dynamic programming techniques that avoid looping over all possible states.  

We hope to investigate ways in which these and other partially observable hybrid systems may be solved using our sufficient statistic, via practical solution strategies.  While such approximate results would still be suboptimal, they may be more informative and accurate than a suboptimal controller that results from a separation between state estimation and control, or from using the belief state given the observations as an (in)sufficient statistic.  We intend to explore approximate solution methods using our sufficient statistic, and compare them to other suboptimal schemes, in future work.  

\begin{ack}
The authors would like to thank the reviewers for their thoroughness and insight, and also David Weirich for his help with functional analysis. 
\end{ack}

\section{Appendix}

\subsection{Proof of Lemma \ref{iidmeas}}
\begin{pf}
The proof follows that of \cite{Elliot93}.
\begin{align*}
\mathbb{P}^{\dag}\left[y_k^x \in A,\, y_k^q = q\mid \mathcal{G}_k\right] &= \mathbb{E}^{\dag}\left[{\bf 1}_A(y_k^x){\bf 1}_{\left\{q\right\}}(y_k^q)\mid \mathcal{G}_k\right] \\
&= \frac{\mathbb{E}^{\dag}\left[{\bf 1}_A(y_k^x){\bf 1}_{\left\{q\right\}}(y_k^q) \Lambda_k^{-1} \mid \mathcal{G}_k\right] }{\mathbb{E}^{\dag}\left[\Lambda_k^{-1} \mid \mathcal{G}_k\right]} \label{fraceq}
\end{align*}
Pulling $\Lambda_{k-1}^{-1}$ outside the expected value from both the numerator and denominator and canceling, since $\Lambda_{k-1}$ is $\mathcal{G}_k$ measurable, the numerator reduces to
\begin{align*}
&\int_A  \frac{\phi(y)\frac{1}{N_q}}{\phi(y - h(x_k,u_{k-1}))Q_{q_k,q}(u_{k-1})}\\
&\hspace{20 mm} \times \mathbb{P}^{\dag} \left[y_k^x = y, y_k^q = q \vert s_k, u_{k-1}\right] dy\\
&=\frac{1}{N_q} \int_A \phi(y) dy = \mathbb{P}^{\dag}\left[y_k^q = q\right]\mathbb{P}^{\dag}\left[y_k^x \in A \right]
\end{align*}
and the denominator becomes
\begin{align*}
&\sum_{q=1}^{N_q}\int_{\mathbb{R}^n} \frac{\phi(y)\frac{1}{N_q}}{\phi(y - h(x_k,u_{k-1}))Q_{q_k,q}(u_{k-1})}\\
& \hspace{20 mm} \times \mathbb{P}^{\dag} \left[y_k^x = y, y_k^q = q \vert s_k, u_{k-1}\right] dy\\
&=\sum_{q=1}^{N_q}\frac{1}{N_q} \int_{\mathbb{R}^n} \phi(y) dy = 1
\end{align*}
Hence,
\begin{equation*}
\mathbb{P}^{\dag}\left[y_k^x \in A,\, y_k^q = q\mid \mathcal{G}_k\right] = \mathbb{P}^{\dag}\left[y_k^q = q\right]\mathbb{P}^{\dag}\left[y_k^x \in A \right]
\end{equation*}
\QEDA
\end{pf}

\subsection{Proof of Lemma \ref{linop}}
\begin{pf}
We first show that $T$ is a bounded linear operator mapping $L^1$ into itself.  We then show that $\sigma_k$ can be defined recursively using $T$.
Linearity follows obviously from the properties of integrals.  For any function $\nu\in L^1(\mathcal{S})$, $u\in\mathcal{U}$, $y\in\mathcal{Y}$,
\begin{align*}
\|T_{u,y}[\nu]\|_{L^1} &= \sum_{q\in\mathcal{Q}}\int_{\mathbb{R}^n}\left|\frac{\phi(y^x - h(x,u))}{\phi(y^x)}N_qQ_{q,y^q}(u)\right.\notag \\
\times &\left. \left[\sum_{q^{-}\in\mathcal{Q}}\int_{\mathbb{R}^n}  {\bf 1}_K(x^-,q^-) \tau(x,q\mid x^-,q^-,u)\right.\right. \notag\\
\end{align*}
\begin{align}
&\hspace{29 mm} \left.\left.\vphantom{\sum_{q^{-}\in\mathcal{Q}}\int_{\mathbb{R}^n}} \times\, \nu(x^-,q^-)\,dx^-\right]\right|dx \notag \\
\|T_{u,y}[\nu]\|_{L^1} &\leq \sum_{q\in\mathcal{Q}} \frac{N_qQ_{q,y^q}(u)}{|\phi(y^x)|}\int_{\mathbb{R}^n}|\phi(y^x - h(x,u))|\,dx\notag \\
&\hspace{10 mm}\times \left[\sum_{q^{-}\in\mathcal{Q}}\int_{\mathbb{R}^n} |\nu(x^-,q^-)|\,dx^-\right] \label{line1}\\
\leq& \sum_{q\in\mathcal{Q}} \frac{N_qQ_{q,y^q}(u)}{|\phi(y^x)|}\sum_{q^{-}\in\mathcal{Q}}\int_{\mathbb{R}^n} |\nu(x^-,q^-)|\,dx^-\label{line2}\\
=&M\|\nu\|_{L^1} \notag
\end{align}
Equation \eqref{line1} follows because ${\bf 1}_K(x^-,q^-) \tau(x,q\mid x^-,q^-,u)\leq 1$ for all $x$, $q$, $x^-$, $q^-$, and \eqref{line2} follows because $\phi$ is a distribution, and $\int\phi(x)\,dx=1$, therefore $\int \phi(y-h(x)) \,dx \leq 1$ for $h$ a bounded continuous function.  Hence for any $\nu\in L^1(\mathcal{S})$, $T$ is a bounded linear operator, with $T_{u,y}[\nu]\in L^1(\mathcal{S})$. 

Induction shows that $\sigma_k=T_{u_{k-1},y_k}[\sigma_{k-1}]$.  Given $\sigma_0 = \rho$, 
\begin{align*}
&T_{u_0,y_1}[\rho] (x,q) = \sum_{q_0\in \mathcal{Q}} \int_{\mathbb{R}^n}{\bf 1}_K(x_0,q_0)\frac{\phi(y_1^x - h(x,u_0))}{\phi(y_1^x)} \\
&\hspace{10 mm}\times  N_qQ_{q,y_1^q}(u_0)\tau(x,q\mid x_0,q_0,u_0) \rho(x_0,q_0)\,dx_0 \\
&=\mathbb{E}^{\dag}\left[{\bf 1}_q(q_1){\bf 1}_x(x_1){\bf 1}_K(x_0,q_0) \Lambda_1 \mid i_1\right] \\
&=\sigma_1(x,q)
\end{align*}
Given $\sigma_l=T_{u_{l-1},y_l}[\sigma_{l-1}] \,\forall \,l=1,\dotsc,k$,
\begin{align*}
&T_{u_{k},y_{k+1}}[\sigma_k] (x,q) = \sum_{q_k\in \mathcal{Q}} \int_{\mathbb{R}^n}{\bf 1}_K(x_k,q_k)\\
 &\hspace{20 mm}\times \frac{\phi(y_{k+1}^x - h(x,u_k))}{\phi(y_{k+1}^x)} N_qQ_{q,y_{k+1}^q}(u_k)\\
&\hspace{25 mm}\times \tau(x,q\mid x_k,q_k,u_k) \sigma_k(x_k,q_k)\,dx_k \\
 &=  \sum_{q_k\in \mathcal{Q}} \int_{\mathbb{R}^n}{\bf 1}_K(x_k,q_k)\frac{\phi(y_{k+1}^x - h(x,u_k))}{\phi(y_{k+1}^x)} \\
&\hspace{20 mm}\times  N_qQ_{q,y_{k+1}^q}(u_k)\tau(x,q\mid x_k,q_k,u_k) \\
 &\hspace{20 mm}\times \left[ \sum_{q_0,\dotsc,q_{k-1}}\int_{\mathbb{R}^n\times\dotsc\times\mathbb{R}^n} \prod_{i=1}^{k-1}{\bf 1}_K(x_i,q_i) \right. \\
&\hspace{20 mm}\times  \tau(x_i,q_i\mid x_{i-1},q_{i-1},u_{i-1}) {\bf 1}_K(x_0,q_0) \\
&\hspace{20 mm}\times  \left.\vphantom{ \sum_{q_0,\dotsc,q_{k-1}}}\rho(x_0,q_0) \Lambda_{k-1}\,dx_0,\dotsc,dx_{k-1}\right] dx_k\\
 \end{align*}
 \begin{align*} 
&=\mathbb{E}^{\dag}\left[{\bf 1}_q(q_{k+1}){\bf 1}_x(x_{k+1})\prod_{i=1}^{k}{\bf 1}_K(x_i,q_i) \right.\\
&\hspace{20 mm} \left.\left. \vphantom{\prod_{i=0}^{k}} \times \Lambda_{k+1}\right|i_{k+1}\right] \\
&=\sigma_{k+1}(x,q)
\end{align*}\QEDA
\end{pf}

\subsection{Proof of Theorem \ref{DPthm}}
The following proofs are based on those appearing in \cite{bertsekas}, chapters 6 and 11.  To facilitate the connection between these proofs and those appearing in \cite{bertsekas} we first reformulate the recursion \eqref{optDP} as a minimization

\begin{equation*}
\sup_{\overline{\pi}} V_0^{\overline{\pi}}(\sigma_0) = -\inf_{\overline{\pi}}- V_0^{\overline{\pi}} = -\inf_{\overline{\pi}}J_0^{\overline{\pi}} 
\end{equation*}

Let $J_k^{\overline{\pi}}(\sigma) = -V_k^{\overline{\pi}}(\sigma)$ and $\overline{\Pi} = \{\overline{\pi} = (\overline{\mu}_0,\,\overline{\mu}_1,\dotsc) : \overline{\mu}_i(\sigma_i) \,\in\, \mathcal{U} \,\forall\, i\}$.  In the following we drop the bar notation over $\overline{\pi}$, and use $\pi$ to denote a policy with respect to the sufficient statistic.  The recursion for $J_k^*(\sigma)$ is identical to that of  $V_k^*(\sigma)$ in \eqref{optDP} except that $J_N^*(\sigma) = -\langle \sigma, {\bf 1}_K \rangle$.

Next we define the operators 
\begin{align*}
H_{\mu}[J] &= \mathbb{E}^{\dag}\left[ J(T_{\mu(\sigma),y}[\sigma])\right] \\
&= \sum_{y^q}\int_{\mathbb{R}^n} J(T_{\mu(\sigma),y}[\sigma])\phi(y^x)\frac{1}{N_q}\,dy^x\\
H[J] &= \inf_{\mu(\sigma)\in\mathcal{U}} H_{\mu}[J]
\end{align*}
The operator $H_{\mu}[J] $ preserves the linearity and boundedness of value function $J(\sigma)$ for all $\sigma$ in $L^1$, which can be seen from a straightforward induction argument.  Because $J_N^{\pi}(\sigma)$ is a bounded linear functional, this then implies that $J_k^{\pi}(\sigma)$ is a bounded linear functional for all $k=0,\dotsc,N$ and for all $\sigma$ in $L^1$.

\begin{lem}\label{ordering}
For all bounded linear functionals $J$, $\sigma \in L^1$, $\pi\in \Pi$, and $r\in\mathbb{R}^+$
\begin{equation*}
H_{\mu}[J] \leq H_{\mu}[J+r] \leq H_{\mu}[J] + r
\end{equation*}
\end{lem}
\begin{pf} 
Because $J \leq J+r$ when $r \geq 0$, we get the following:
\begin{align*}
H_{\mu}[J](y) &= \sum_{y^q} \int_{\mathbb{R}^n} J(T_{\mu(\sigma),y} [\sigma])\frac{1}{N_q}\phi(y^x)dy^x \\
&\leq \sum_{y^q} \int_{\mathbb{R}^n} (J(T_{\mu(\sigma),y} [\sigma])+r)\frac{1}{N_q}\phi(y^x)dy^x\\
&=H_{\mu}[J+r](y) \\
&=\sum_{y^q} \int_{\mathbb{R}^n} J(T_{\mu(\sigma),y} [\sigma])\frac{1}{N_q}\phi(y^x)dy^x + r \\
&=H_{\mu}[J]+r
\end{align*}\QEDA
\end{pf}

\begin{prop}\label{contractprop}
For any $M \in \mathbb{N}$, where $J_0^*(\sigma) = \inf_{\pi \in \Pi }J_0^{\pi}(\sigma)$,
\begin{equation*}
J_0^*(\sigma) = H^M[J_{M}^*(\sigma)]
\end{equation*}
 Further, for any $\epsilon >0$ there exists an M-stage $\epsilon$-optimal policy $\pi_{\epsilon}$, defined as
\begin{equation*}
J_0^* \leq J_0^{\pi_{\epsilon}} \leq J_0^* + \epsilon
\end{equation*}
\end{prop}

\begin{pf}{(By backwards induction on $M$)}

For $M=N$, 
\begin{equation*}
J_{N}^* (\sigma) = H^0[J_{N}^*(\sigma)] \end{equation*}
(since $H^0[J] = J$).  Also, because $J_{N}$ does not depend on a control input, $J_{N}^*(\sigma) = J_{N}^{\pi}(\sigma)$ for any policy $\pi\in \Pi$.  Therefore, for any $\epsilon>0$, $J_{N}^{\pi_{\epsilon}}(\sigma) = J_{N}^*(\sigma) \leq J_{N}^*(\sigma) + \epsilon$.

Assume for $M=k+1$ that $J_{k+1}^*(\sigma)= H^{N-k-1}[J_{N}^*(\sigma)]$ and that for all $\epsilon > 0 $ there exists a policy $\pi_{\epsilon}$ such that $J_{k+1}^{\pi_{\epsilon}}(\sigma) \leq J_{k+1}^*(\sigma) + \epsilon$.
Then by Lemma \ref{ordering}, for any $\mu \in \Pi$
\begin{equation*}
H_{\mu}[J_{k+1}^{\pi_{\epsilon}}] \leq H_{\mu}[J_{k+1}^* + \epsilon] \leq H_{\mu}[J_{k+1}^*] + \epsilon \end{equation*}
By aggregating $\mu$ with the control function $\pi_{\epsilon}$ to get $\hat{\pi}_{\epsilon} = (\mu, \pi_{\epsilon} )$, we then have
\begin{align*}
\inf_{\pi} J_{k}^{\pi}(\sigma) \leq J_{k}^{\hat{\pi}_{\epsilon} }(\sigma) &= H_{\mu}[J_{k+1}^{\pi_{\epsilon}}] \\
	&\leq H_{\mu}[J_k^*] + \epsilon \\
	\end{align*}
Since the above holds for  \emph{any} $\mu \in \Pi$, 
\begin{align*}
\inf_{\pi} J_{k}^{\pi}(\sigma) &\leq H[J_{k+1}^*(\sigma)] \\
	&= H[H^{N-k-1}[J_{N}^*(\sigma)]] = H^{N-k}[J_{N}^*(\sigma)]
\end{align*}
By definition $H^{N-k}[J_{N}^*(\sigma)] \leq J_{k}^*(\sigma)$, hence \newline$H^{N-k}[J_{N}^*(\sigma)] =J_{k}^*(\sigma)$.

Next, by the induction argument, for any $\hat{\epsilon}>0$, let $\hat{\pi}$ be such that 
\begin{equation*}
J_{k+1}^{\hat{\pi}}(\sigma) \leq J_{k+1}^*(\sigma) + \frac{\hat{\epsilon}}{2}	\end{equation*}
Let $\hat{\mu} \in \Pi$ be such that
\begin{equation*}
H_{\hat{\mu}}[J_{k+1}^*] \leq H[J_{k+1}^*] + \frac{\overline{\epsilon}}{2} 
\end{equation*}
Define $\hat{\pi}_{\epsilon} = (\hat{\mu},\hat{\pi})$.  Then
\begin{align*}
J_{k}^{\hat{\pi}_{\epsilon}} =H_{\hat{\mu}}[J_{k+1}^{\hat{\pi}}] &\leq H_{\hat{\mu}}[J_{k+1}^*] + \frac{\hat{\epsilon}}{2} \\
&\leq H[J_{k+1}^*] + \frac{\hat{\epsilon}}{2} + \frac{\hat{\epsilon}}{2} \\
&= J_{k}^* + \hat{\epsilon}
\end{align*}
It follows from induction that $J_0^* \leq J_0^{\pi_{\epsilon}} \leq J_0^* + \epsilon$ for any $M$.\QEDA
\end{pf}
We also use the result from \cite{bertsekas} on the existence of a uniformly N-stage optimal policy $\pi^* = (\mu_0^*,\mu_1^*,...)$, which we give without proof (see Ch. 6), since the proof does not change in our context.
\begin{prop}\label{noptprop}
A policy is uniformly N-stage optimal if and only if $H_{\mu_k^*}[H^{N-k-1}[J_{N}^*]] = H^{N-k}[J_{N}^*]$ for all $k=0,\dotsc,N$, and this policy exists if and only if the infimum of
\begin{equation*}
H^{N-k}[J_{N}^*] = \inf_{u\in\mathcal{U}}H_u\left[H^{N-k-1}[J_N^*]\right]	\end{equation*}
is attained for all $\sigma \in L^1$ and $k=0,\dotsc,N$.  A sufficient condition for the infimum to be attained is that
\begin{equation*}
U_k(\sigma,\lambda) = \left\{ u \in \mathcal{U} : H_u\left[H^{N-k-1}[J_N^*]\right] \leq \lambda \right\}
\end{equation*}
is compact for all $\sigma\in L^1,\, \lambda \in \mathbb{R}$, and $k=0,\dotsc,N$.
\end{prop}

Substituting $V_k^*= -J_k^*$, it is clear that Prop. \ref{contractprop} validates the dynamic programming algorithm \eqref{optDP}, and proves the existence of at least an $\epsilon$-optimal policy, and so the first part of Theorem \ref{DPthm} is proved.  Finally, using Prop. \ref{noptprop}, because $\mathcal{U}$ is defined as a compact (i.e. closed and bounded) Borel set, and $J_k$ (and so $V_k$) is bounded for all $\sigma \in L^1$ and for each $u\in\mathcal{U}$, then there exists some $u\in\mathcal{U}$ such that the infimum in $\inf_{u\in\mathcal{U}}H_u\left[H^{N-k-1}[J_N^*]\right]$ $=\inf_{u\in\mathcal{U}}\mathbb{E}[J_k(T_{u,y}[\sigma])]$ is attained for all $k$ (and likewise the supremum of $V_k$ is achieved for all $k$).  Therefore, for \eqref{optDP}, there always exists an optimal policy $\pi$ given by \eqref{DPpol}.

\subsection{Proof of Theorem \ref{equival}}
\begin{pf}
For a vector $\overline{u}=[u_0, u_1,\dotsc u_{N-1}]$ with each $u_i\in\mathcal{U}$, we have by definition that 
\begin{equation*}
r_K(\overline{u}) = \overline{r}_K(\overline{u}) \hspace{5 mm}\forall\, \overline{u}\in\mathcal{U}_N
\end{equation*}
 Since $\sigma_k = \eta_k(\rho, i_k)$, the control policy $\overline{\pi}=\newline (\overline{\mu}_0(\sigma_0),\overline{\mu}_1(\sigma_1),\dotsc)$ can be rewritten as a function of the information vector $i_k$, where $\overline{\mu}_k(\sigma_k) = \overline{\mu}_k(\eta_k(\rho, i_k)) = \mu_k(i_k)$.  Then by defining the policy $\pi$ in terms of $\mu$, we have that $r_K(\pi) = \overline{r}_K(\overline{\pi})$ for all $\overline{\pi}\in\overline{\Pi}$.  If $\overline{\pi}^*$ is optimal for $\overline{r}_K(\overline{\pi})$, it then must be optimal for $r_K(\pi)$ as well, and further,
 \begin{equation*}
 r_K(\pi^*) = \overline{r}_K(\overline{\pi}^*)
 \end{equation*}\QEDA
 \end{pf}

\bibliographystyle{elsarticle-harv}
\bibliography{hybridsystems}

\begin{thebibliography}{31}
\expandafter\ifx\csname natexlab\endcsname\relax\def\natexlab#1{#1}\fi
\providecommand{\url}[1]{\texttt{#1}}
\providecommand{\href}[2]{#2}
\providecommand{\path}[1]{#1}
\providecommand{\DOIprefix}{doi:}
\providecommand{\ArXivprefix}{arXiv:}
\providecommand{\URLprefix}{URL: }
\providecommand{\Pubmedprefix}{pmid:}
\providecommand{\doi}[1]{\href{http://dx.doi.org/#1}{\path{#1}}}
\providecommand{\Pubmed}[1]{\href{pmid:#1}{\path{#1}}}
\providecommand{\bibinfo}[2]{#2}
\ifx\xfnm\relax \def\xfnm[#1]{\unskip,\space#1}\fi
\bibitem[{Abate et~al.(2007)Abate, Amin, Prandini, Lygeros and
  Sastry}]{Abatecomp}
\bibinfo{author}{Abate, A.}, \bibinfo{author}{Amin, S.},
  \bibinfo{author}{Prandini, M.}, \bibinfo{author}{Lygeros, J.},
  \bibinfo{author}{Sastry, S.}, \bibinfo{year}{2007}.
\newblock \bibinfo{title}{Computational approches to reachability analysis of
  stochastic hybrid systems}, in: \bibinfo{booktitle}{Hybrid Systems:
  Computation and Control}.
\bibitem[{Abate et~al.(2008)Abate, Prandini, Lygeros and Sastry}]{Abate1}
\bibinfo{author}{Abate, A.}, \bibinfo{author}{Prandini, M.},
  \bibinfo{author}{Lygeros, J.}, \bibinfo{author}{Sastry, S.},
  \bibinfo{year}{2008}.
\newblock \bibinfo{title}{Probabilistic reachability and safety for controlled
  discrete time stochastic hybrid systems}.
\newblock \bibinfo{journal}{Automatica} \bibinfo{volume}{44},
  \bibinfo{pages}{2724--2734}.
\bibitem[{Aoki(1989)}]{Aoki}
\bibinfo{author}{Aoki, M.}, \bibinfo{year}{1989}.
\newblock \bibinfo{title}{Optimization of Stochastic Systems: Topics in
  Discrete-Time Dynamics}.
\newblock \bibinfo{publisher}{Academic Press, Inc.}
\bibitem[{Bai et~al.(2011)Bai, Hsu, Lee and Ngo}]{Bai2011}
\bibinfo{author}{Bai, H.}, \bibinfo{author}{Hsu, D.}, \bibinfo{author}{Lee,
  W.}, \bibinfo{author}{Ngo, V.A.}, \bibinfo{year}{2011}.
\newblock \bibinfo{title}{Monte carlo value iteration for continuous-state
  {POMDP}s}, in: \bibinfo{booktitle}{Algorithmic Foundations of Robotics IX}.
  volume~\bibinfo{volume}{68}, pp. \bibinfo{pages}{175--191}.
\bibitem[{Bertsekas and Shreve(1996)}]{bertsekas}
\bibinfo{author}{Bertsekas, D.P.}, \bibinfo{author}{Shreve, S.E.},
  \bibinfo{year}{1996}.
\newblock \bibinfo{title}{Stochastic Optimal Control: The Discrete-Time Case}.
\newblock \bibinfo{publisher}{Athena Scientific}.
\bibitem[{Brunskill et~al.(2008)Brunskill, Kaelbling, Lozano-Perez and
  Roy}]{Roy}
\bibinfo{author}{Brunskill, E.}, \bibinfo{author}{Kaelbling, L.},
  \bibinfo{author}{Lozano-Perez, T.}, \bibinfo{author}{Roy, N.},
  \bibinfo{year}{2008}.
\newblock \bibinfo{title}{Continuous-state {P}{O}{M}{D}{P}s with hybrid
  dynamics}, in: \bibinfo{booktitle}{Symposium on Artificial Inteligence and
  Mathematics}.
\bibitem[{Bujorianu(2010)}]{Bujor2}
\bibinfo{author}{Bujorianu, M.L.}, \bibinfo{year}{2010}.
\newblock \bibinfo{title}{Variational inequalities for the stochastic
  reachability problem}, in: \bibinfo{booktitle}{{IEEE} Conference on Decision
  and Control}, pp. \bibinfo{pages}{1854--1859}.
\bibitem[{Caldwell and Murphy(2011)}]{caldwell}
\bibinfo{author}{Caldwell, T.}, \bibinfo{author}{Murphy, T.},
  \bibinfo{year}{2011}.
\newblock \bibinfo{title}{Switching mode generation and optimal estimation with
  application to skid-steering}.
\newblock \bibinfo{journal}{Automatica} \bibinfo{volume}{47},
  \bibinfo{pages}{50--64}.
\bibitem[{Ding et~al.(2013)Ding, Abate and Tomlin}]{Ding2013}
\bibinfo{author}{Ding, J.}, \bibinfo{author}{Abate, A.},
  \bibinfo{author}{Tomlin, C.}, \bibinfo{year}{2013}.
\newblock \bibinfo{title}{Optimal control of partially observable discrete time
  stochastic hybrid systems for safety specifications}, in:
  \bibinfo{booktitle}{American Control Conference}, pp.
  \bibinfo{pages}{6231--6236}.
\bibitem[{Elliot(1993)}]{Elliot93}
\bibinfo{author}{Elliot, R.J.}, \bibinfo{year}{1993}.
\newblock \bibinfo{title}{A general recursive discrete-time filter}.
\newblock \bibinfo{journal}{Journal of Applied Probability}
  \bibinfo{volume}{30}, \bibinfo{pages}{575--588}.
\bibitem[{Esfahani et~al.(2011)Esfahani, Chaterjee and Lygeros}]{Esfahani}
\bibinfo{author}{Esfahani, P.M.}, \bibinfo{author}{Chaterjee, D.},
  \bibinfo{author}{Lygeros, J.}, \bibinfo{year}{2011}.
\newblock \bibinfo{title}{On a problem of stochastic reach-avoid set
  characterization}, in: \bibinfo{booktitle}{{IEEE} Conference on Decision and
  Control and European Control Conference}, pp. \bibinfo{pages}{7069--7074}.
\bibitem[{Fernandez-Gaucherand and Marcus(1997)}]{Gauch}
\bibinfo{author}{Fernandez-Gaucherand, E.}, \bibinfo{author}{Marcus, S.I.},
  \bibinfo{year}{1997}.
\newblock \bibinfo{title}{Risk-sensitive optimal control of hidden markov
  models: Structural results}.
\newblock \bibinfo{journal}{{IEEE} Transactions on Automatic Control}
  \bibinfo{volume}{42}, \bibinfo{pages}{1418--1422}.
\bibitem[{Hofbaur and Williams(2004)}]{hof1}
\bibinfo{author}{Hofbaur, M.W.}, \bibinfo{author}{Williams, B.C.},
  \bibinfo{year}{2004}.
\newblock \bibinfo{title}{Hybrid estimation of complex systems}.
\newblock \bibinfo{journal}{{IEEE} Transactions on Systems, Man, and
  Cybernetics} \bibinfo{volume}{34}, \bibinfo{pages}{2178--2191}.
\bibitem[{Hu et~al.(2000)Hu, Lygeros and Sastry}]{Hu1}
\bibinfo{author}{Hu, J.}, \bibinfo{author}{Lygeros, J.},
  \bibinfo{author}{Sastry, S.}, \bibinfo{year}{2000}.
\newblock \bibinfo{title}{Towards a theory of stochastic hybrid systems}.
\newblock \bibinfo{journal}{Hybrid Systems: Computation and Control} ,
  \bibinfo{pages}{160--173}.
\bibitem[{James et~al.(1994)James, Baras and Elliot}]{James}
\bibinfo{author}{James, M.R.}, \bibinfo{author}{Baras, J.S.},
  \bibinfo{author}{Elliot, R.J.}, \bibinfo{year}{1994}.
\newblock \bibinfo{title}{Risk-sensitive control and dynamic games for
  partially observed discrete-time nonlinear systems}.
\newblock \bibinfo{journal}{{IEEE} Transactions on Automatic Control}
  \bibinfo{volume}{39}, \bibinfo{pages}{780--792}.
\bibitem[{Kamgarpour et~al.(2011)Kamgarpour, Ding, Summers, Abate, Lygeros and
  Tomlin}]{Kamgar1}
\bibinfo{author}{Kamgarpour, M.}, \bibinfo{author}{Ding, J.},
  \bibinfo{author}{Summers, S.}, \bibinfo{author}{Abate, A.},
  \bibinfo{author}{Lygeros, J.}, \bibinfo{author}{Tomlin, C.},
  \bibinfo{year}{2011}.
\newblock \bibinfo{title}{Discrete time stochastic hybrid dynamical games:
  Verification and controller synthesis}, in: \bibinfo{booktitle}{{IEEE}
  Conference on Decision and Control}.
\bibitem[{Kariotoglou et~al.(2013)Kariotoglou, Summers, Summers, Kamgarpour and
  Lygers}]{Kariot2013}
\bibinfo{author}{Kariotoglou, N.}, \bibinfo{author}{Summers, S.},
  \bibinfo{author}{Summers, T.}, \bibinfo{author}{Kamgarpour, M.},
  \bibinfo{author}{Lygers, J.}, \bibinfo{year}{2013}.
\newblock \bibinfo{title}{Approximate dynamic programming for stochastic
  reachability}, in: \bibinfo{booktitle}{European Control Conference}, pp.
  \bibinfo{pages}{584 -- 589}.
\bibitem[{Koutsoukos et~al.(2003)Koutsoukos, Kurien and Zhao}]{Kouts1}
\bibinfo{author}{Koutsoukos, X.}, \bibinfo{author}{Kurien, J.},
  \bibinfo{author}{Zhao, F.}, \bibinfo{year}{2003}.
\newblock \bibinfo{title}{Estimation of distributed hybrid systems using
  particle filtering methods}, in: \bibinfo{booktitle}{Hybrid Systems:
  Computation and Control}, pp. \bibinfo{pages}{298--313}.
\bibitem[{Liu and Hwang(2012)}]{Hwang1}
\bibinfo{author}{Liu, W.}, \bibinfo{author}{Hwang, I.}, \bibinfo{year}{2012}.
\newblock \bibinfo{title}{A stochastic approximation based state estimation
  algorithm for stochastic hybrid systems}, in: \bibinfo{booktitle}{American
  Control Conference}, pp. \bibinfo{pages}{312--317}.
\bibitem[{Mitchell(2008)}]{Mitchell2}
\bibinfo{author}{Mitchell, I.}, \bibinfo{year}{2008}.
\newblock \bibinfo{title}{The flexible, extensible and efficient toolbox of
  level set methods}.
\newblock \bibinfo{journal}{Journal of Scientific Computing}
  \bibinfo{volume}{35}, \bibinfo{pages}{300--329}.
\newblock \bibinfo{note}{10.1007/s10915-007-9174-4}.
\bibitem[{Porta et~al.(2006)Porta, Vlassis, Spain and Poupart}]{Porta}
\bibinfo{author}{Porta, J.M.}, \bibinfo{author}{Vlassis, N.},
  \bibinfo{author}{Spain, M.T.}, \bibinfo{author}{Poupart, P.},
  \bibinfo{year}{2006}.
\newblock \bibinfo{title}{Point-based value iteration for continuous
  {P}{O}{M}{D}{P}s}.
\newblock \bibinfo{journal}{Journal of Machine Learning Research}
  \bibinfo{volume}{7}, \bibinfo{pages}{2329--2367}.
\bibitem[{Prandini and Hu(2006)}]{Prand1}
\bibinfo{author}{Prandini, M.}, \bibinfo{author}{Hu, J.}, \bibinfo{year}{2006}.
\newblock \bibinfo{title}{Stochastic Reachability: Theoretical Foundations and
  Numerical Approximation}. \bibinfo{publisher}{Springer Verlag}.
\newblock Lecture Notes in Control and Information Sciences, pp.
  \bibinfo{pages}{107--139}.
\bibitem[{Shiryaev(1964)}]{shiryaev}
\bibinfo{author}{Shiryaev, A.}, \bibinfo{year}{1964}.
\newblock \bibinfo{title}{On markov sufficient statistics in nonadditive
  {B}ayes problems of sequential analysis}.
\newblock \bibinfo{journal}{Theory of Probability and its Applications}
  \bibinfo{volume}{9}, \bibinfo{pages}{604--618}.
\bibitem[{Soudjani and Abate(2013)}]{SoudjaniSiam13}
\bibinfo{author}{Soudjani, S.}, \bibinfo{author}{Abate, A.},
  \bibinfo{year}{2013}.
\newblock \bibinfo{title}{Adaptive and sequential gridding procedures for the
  abstraction and verification of stochastic processes}.
\newblock \bibinfo{journal}{{SIAM} Journal on Applied Dynamical Systems}
  \bibinfo{volume}{12}, \bibinfo{pages}{921--956}.
\bibitem[{Stein and Shakarchi(2005)}]{SteinShak}
\bibinfo{author}{Stein, E.}, \bibinfo{author}{Shakarchi, R.},
  \bibinfo{year}{2005}.
\newblock \bibinfo{title}{Real Analysis: Measure Theory, Integration, and
  Hilbert Spaces}.
\newblock Princeton Lectures in Analysis, \bibinfo{publisher}{Princeton
  University Press}.
\bibitem[{Summers and Lygeros(2010)}]{summers}
\bibinfo{author}{Summers, S.}, \bibinfo{author}{Lygeros, J.},
  \bibinfo{year}{2010}.
\newblock \bibinfo{title}{Verification of discrete time stochastic hybrid
  systems: A stochastic reach-avoid decision problem}.
\newblock \bibinfo{journal}{Automatica} \bibinfo{volume}{46},
  \bibinfo{pages}{1951--1961}.
\bibitem[{Tkachev et~al.(2013)Tkachev, Katoen, Mereacre and
  Abate}]{Tkachev2013}
\bibinfo{author}{Tkachev, I.}, \bibinfo{author}{Katoen, J.P.},
  \bibinfo{author}{Mereacre, A.}, \bibinfo{author}{Abate, A.},
  \bibinfo{year}{2013}.
\newblock \bibinfo{title}{Quantitative automata-based controller synthesis for
  non-autonomous stochastic hybrid systems}, in: \bibinfo{booktitle}{Hybrid
  Systems: Computation and Control}, pp. \bibinfo{pages}{293--302}.
\bibitem[{Tomlin et~al.(2003)Tomlin, Mitchell, Bayen and Oishi}]{oishi1}
\bibinfo{author}{Tomlin, C.}, \bibinfo{author}{Mitchell, I.},
  \bibinfo{author}{Bayen, A.}, \bibinfo{author}{Oishi, M.},
  \bibinfo{year}{2003}.
\newblock \bibinfo{title}{Computational techniques for the verification and
  control of hybrid systems}, in: \bibinfo{booktitle}{Proceedings of the
  {IEEE}}, pp. \bibinfo{pages}{986--1001}.
\bibitem[{Verma and del Vecchio(2009)}]{Verm1}
\bibinfo{author}{Verma, R.}, \bibinfo{author}{del Vecchio, D.},
  \bibinfo{year}{2009}.
\newblock \bibinfo{title}{Continuous control of hybrid automata with imperfect
  mode information assuming separation between state estimation and control},
  in: \bibinfo{booktitle}{{IEEE} Conference on Decision and Control},
  \bibinfo{address}{Shanghai}.
\bibitem[{Verma and del Vecchio(2010)}]{Verm2}
\bibinfo{author}{Verma, R.}, \bibinfo{author}{del Vecchio, D.},
  \bibinfo{year}{2010}.
\newblock \bibinfo{title}{Control of hybrid automata with hidden modes:
  Translation to a perfect state information problem}, in:
  \bibinfo{booktitle}{{IEEE} Conference on Decision and Control}.
\bibitem[{Vrakopoulou et~al.(2013)Vrakopoulou, Margellos, Lygeros and
  Andersson}]{Vrakop2013}
\bibinfo{author}{Vrakopoulou, M.}, \bibinfo{author}{Margellos, K.},
  \bibinfo{author}{Lygeros, J.}, \bibinfo{author}{Andersson, G.},
  \bibinfo{year}{2013}.
\newblock \bibinfo{title}{Probabilistic guarantees for the n-1 security of
  systems with wind power generation}, in: \bibinfo{booktitle}{Reliability and
  Risk Evaluation of Wind Integrated Power Systems}.
  \bibinfo{publisher}{Springer India}, pp. \bibinfo{pages}{59--73}.

\end{thebibliography}

\end{document}